\documentclass{amsart}
\usepackage{graphicx}
\usepackage{amsfonts}
\newtheorem{tm}{Theorem}

\newtheorem{defi}{Definition}

\newtheorem{rem}{Remark}
\newtheorem{rems}{Remarks}

\newtheorem{ex}{Example}

\newtheorem{nota}{Notation}

\begin{document}

\title{Hyperbolic polynomials and rigid moduli orders}
\author{Vladimir Petrov Kostov}
\address{Universit\'e C\^ote d’Azur, CNRS, LJAD, France} 
\email{vladimir.kostov@unice.fr}

\begin{abstract}
  A hyperbolic polynomial (HP) is a real univariate polynomial with all roots
  real. By Descartes' rule of signs a HP with all coefficients nonvanishing
  has exactly $c$ positive and exactly $p$ negative roots counted with
  multiplicity, where $c$ and $p$ are the numbers of sign changes and sign
  preservations in the sequence of its coefficients. We consider HPs with
  distinct moduli of the roots.
  We ask the question when the order of the moduli of the negative roots
  w.r.t. the positive roots on the real positive half-line completely
  determines the signs of the coefficients of the polynomial. When there is
  at least one positive and at least one negative root this is possible
  exactly when the moduli of the negative roots interlace with the
  positive roots (hence half or about half of the roots are positive). In this
  case the signs of the coefficients of the HP are either
  $(+,+,-,-,+,+,-,-,\ldots )$ or $(+,-,-,+,+,-,-,+,\ldots )$.\\ 

  {\bf Key words:} real polynomial in one variable; hyperbolic polynomial; sign
  pattern; Descartes' 
rule of signs\\ 

{\bf AMS classification:} 26C10; 30C15
\end{abstract}
\maketitle

\section{Introduction}

We consider {\em hyperbolic polynomials (HPs)}, i.e. real univariate
polynomials with all roots real. We assume the leading coefficient
to be positive and all coefficients to be nonvanishing. The Descartes' rule
of signs applied to such a degree $d$ HP with $c$ sign changes and $p$
sign preservations
in the sequence of its coefficients, $c+p=d$, implies that
the HP has $c$ positive and
$p$ negative roots counted with multiplicity. In what follows we consider
the {\em generic} case when the moduli of all roots are distinct.

\begin{defi}\label{defiSPMO}
  {\rm (1) A {\em sign pattern (SP)} of length $d+1$
    is a sequence of $d+1$ $(+)$- and/or $(-)$-signs.
    We say that the polynomial $Q:=x^d+\sum _{j=0}^{d-1}a_jx^j$ defines
    (or realizes) the SP
    $\sigma (Q):=(+,{\rm sgn}(a_{d-1}),\ldots ,{\rm sgn}(a_0))$.

    (2) A {\em moduli order (MO)} of length $d$ is a formal string of $c$
    letters $P$ and $p$ letters $N$ separated by signs of inequality $<$.
    These letters indicate the relative positions of the moduli
    of the roots of the HP on the real positive half-line. 
    E.g. for $d=6$, to say that a given HP $Q$ defines (or realizes)
    the MO
    $N<N<P<N<P<N$
    means that for $\sigma (Q)$, one has $c=2$ and $p=4$ and that for
    the positive roots $\alpha _1<\alpha _2$ and the negative roots
    $-\gamma _j$ of $Q$, one has
    $\gamma _1<\gamma _2<\alpha _1<\gamma _3<\alpha _2<\gamma _4$.

    (3) We say that a given MO {\em realizes} a given SP if there
    exists a HP which defines the given MO and the given~SP.}
  \end{defi}

\begin{ex}\label{ex1}
  {\rm For $d=1$, if the SP defined by a HP with a nonzero root equals
$(+,+)$ (resp. $(+,-)$),  then this root is negative (resp. positive).

    For $d=2$, a HP with roots of opposite signs and different moduli
    defines the SP $(+,+,-)$ with MO $P<N$ or the SP $(+,-,-)$
with MO $N<P$.}
  \end{ex}

\begin{rem}\label{remconcat}
  {\rm Suppose that the MO $r$ is realizable by a HP $Q$.
    Denote by $rP$, $rN$ (resp. $Pr$ and $Nr$)
    the MOs obtained from $r$ by adding to the right the inequality
    $<P$ or $<N$ (resp. by adding to the left the inequality $P<$ or $N<$).
    For $\varepsilon >0$ sufficiently small, the product
    $(x-\varepsilon )Q(x)$ (resp. $(x+\varepsilon )Q(x)$) defines the
    MO $Pr$ (resp. $Nr$). Indeed, the modulus of the root
    $\pm \varepsilon$ is much smaller than any of the moduli of the roots
    of~$Q$. In the same way, the product
    $-(1-\varepsilon x)Q(x)$ (resp. $(1+\varepsilon x)Q(x)$) defines the
    MO $rP$ (resp. $rN$), because the modulus of the
    root $\pm 1/\varepsilon$ is much larger than any of the moduli of the roots
    of~$Q$.  When several products of the form $(x\pm \varepsilon )Q(x)$ and/or $\pm (1\pm \varepsilon x)Q(x)$ are used, then they are performed with different numbers $\varepsilon _j$ for which one has $0<\cdots \ll \varepsilon _{j+1}\ll \varepsilon _j$.
  } 
  \end{rem}

\begin{defi}\label{defirigid}
  {\rm A MO is {\em rigid} if all HPs realizing this MO
    define one and the same SP, i.e. if the MO realizes only one~SP.}
  \end{defi}

The aim of the present paper is to characterize all rigid MOs.
From now on we assume that $c\geq 1$ and $p\geq 1$. Indeed, when all roots
are of the same sign, then there is a single SP corresponding to such a MO 
(this is either the all-pluses SP when the roots are negative or $(+,-,+,-,+,\ldots )$ when they are positive), so according to
our definition this MO should be considered as rigid. However as
it excludes the question how moduli of negative roots are placed w.r.t.
the positive roots on the real positive half-line, this case should be
considered as trivial.

\begin{nota}\label{nota1}
  {\rm (1) We introduce the following four MOs:}

  $$\begin{array}{lcl}
    r_{PN}:P<N<P<N<\cdots <N~,&&r_{PP}:P<N<P<N<\cdots <P~,\\ \\
    
    r_{NP}:N<P<N<P<\cdots <P&{\rm and}&r_{NN}:N<P<N<P<\cdots <N~.
  \end{array}$$
  {\rm The orders $r_{PN}$ and $r_{NP}$ (resp. $r_{PP}$ and $r_{NN}$)
    correspond to even (resp. to odd) degree $d$. In the case of $r_{PN}$ and
    $r_{NP}$ there are $d/2$ positive and $d/2$ negative roots, in the case of 
    $r_{PP}$ there are $(d+1)/2$ positive and $(d-1)/2$ negative roots and
    vice versa in the case of $r_{NN}$. The MOs $r_{PN}$, $r_{NP}$,
    $r_{PP}$ and $r_{NN}$ are the only ones in which there are no
    two consecutive moduli of roots of one and the same sign hence
    for $d\geq 3$, they are the ones and the only ones which contain no (sub)string
    of the form $P<P<N$, $N<N<P$, $N<P<P$ or $P<N<N$. 

    (2) We are particularly
    interested in the following two SPs:}

  $$
  \Sigma _{+}:=(+,+,-,-,+,+,-,-,\ldots )~~~\, {\rm and}~~~\, 
  \Sigma _{-}:=(+,-,-,+,+,-,-,+,\ldots )~.$$
      
  \end{nota}

The main result of the paper is the following theorem:

\begin{tm}\label{tmmain}
  (1) For $d\geq 3$, a MO different from $r_{PN}$, $r_{NP}$,
    $r_{PP}$ and $r_{NN}$ is not rigid.

  (2) For $d\geq 1$, the MOs $r_{PP}$, $r_{PN}$, $r_{NP}$ and $r_{NN}$ are rigid.
  When the roots of a HP define one of these MOs, then the SP
  of the HP is one of the SPs $\Sigma _{\pm}$. The exact correspondence
  is given by the following table (its fourth and seventh
  columns contain the last three
  signs of the SP; the degree $d$ is considered modulo~$4$):

  $$\begin{array}{ccccccccc}
    d\mod(4)&&{\rm MO}&{\rm SP}&&~~~\, \, &{\rm MO}&{\rm SP}&\\ \\
    0&&r_{NP}&\Sigma _-&-~+~+&&r_{PN}&\Sigma _+&-~-~+\\ \\
    1&&r_{PP}&\Sigma _-&+~+~-&&r_{NN}&\Sigma _+&-~+~+\\ \\
    2&&r_{NP}&\Sigma _-&+~-~-&&r_{PN}&\Sigma _+&+~+~-\\ \\
    3&&r_{PP}&\Sigma _-&-~-~+&&r_{NN}&\Sigma _+&+~-~-\end{array}$$
\end{tm}

The theorem is proved in Section~\ref{secprtmmain}. 
Our next step is to consider the possibility to have equalities between
moduli of roots and zeros among the coefficients.

\begin{rem}
  {\rm We remind that a HP $Q$ with nonvanishing constant term
cannot have two consecutive vanishing coefficients. Indeed, if $Q$ is
hyperbolic, then its nonconstant derivatives are also hyperbolic and the
{\em reverted polynomial} $x^dQ(1/x)$ is also hyperbolic. Suppose that
$Q$ is hyperbolic and has two or more consecutive vanishing coefficients.
Then applying derivation
and reversion to $Q$ one can obtain a polynomial of the form $Ax^s+B$,
$s\geq 3$, $A$, $B\in \mathbb{R}^*$, which must be hyperbolic. This, however,
is impossible.}
\end{rem}

\begin{defi}
  {\rm A {\em sign pattern admitting zeros (SPAZ)} of length $d+1$ is a
    sequence of $d+1$ $(+)$- and/or $(-)$-signs and eventually zeros. The first
    element of the sequence must be a $(+)$-sign. To determine the number of
    sign changes and sign  preservations of a SPAZ one has to erase the zeros. 
    A {\em moduli order admitting equalities (MOAE)} of length $d$ is a
    formal string of letters $P$ and $N$ separated by signs of inequality
    $\leq$. E.g. for $d=6$, saying that the HP $Q$ defines the MOAE
    $N\leq N\leq P\leq N\leq P\leq N$ means that the SPAZ defined by $Q$ and the one defined by $Q(-x)$ have 
    at least two and four sign changes respectively, and the constant term of $Q$ is
    nonvanishing; for the moduli of its roots (with the notation from
    Definition~\ref{defiSPMO}), one has
    $\gamma _1\leq \gamma _2\leq \alpha _1\leq
    \gamma _3\leq \alpha _2\leq \gamma _4$.} 
\end{defi}

\begin{ex}\label{ex2}
  {\rm For $d=2$, a HP with nonvanishing constant term and two opposite roots
    is of the form $F:=x^2-a^2$, $a\in \mathbb{R}^*$. It defines the SPAZ
    $(+,0,-)$ which has one sign change and no sign preservation. One has $F(-x)=F(x)$.}
  \end{ex}

\begin{nota}
  {\rm We denote by $r_{PN}^0$, $r_{PP}^0$, $r_{NN}^0$ and $r_{NP}^0$ the MOsAE
    obtained from the respective MOs $r_{PN}$, $r_{PP}$, $r_{NN}$ and $r_{NP}$
    (see Notation~\ref{nota1}) by replacing
    the inequalities $<$ by inequalities~$\leq$.}
\end{nota}

\begin{tm}\label{tm2}
  (1) Suppose that $d\geq 1$ is odd. If a HP with nonvanishing constant term
  defines the MOAE $r_{PP}^0$ or 
  $r_{NN}^0$, then this HP has no vanishing coefficient
  and defines the SP
  as claimed by part (2) of Theorem~\ref{tmmain}.

  (2) If $d\geq 2$ is even and a HP with nonvanishing constant term defines
  the MOAE $r_{PN}^0$ or 
  $r_{NP}^0$, then either

  (i) this HP is even hence of the form $A\prod _{j=1}^{d/2}(x^2-a_j^2)$, where $A>0$, 
  $a_j\in \mathbb{R}^*$ are not necessarily distinct and the HP defines the SPAZ
  $(+,0,-,0,+,0,-,0,\ldots )$, or 

  (ii) this HP has no vanishing coefficient, it defines
  the SP as claimed by part (2) of Theorem~\ref{tmmain} and
  it is not possible to
  represent the set of its roots as a union of $d/2$ couples of the form
  $\{ a_j, -a_j\}$.
    \end{tm}

The theorem is proved in Section~\ref{secprtm2}. In the next section we
compare the problem to characterize rigid MOs to other problems
arising in the theory of real univariate polynomials.

\section{Other related problems}

A rigid MO is one which uniquely defines the SP. One could ask the inverse
question, whether there exist SPs which uniquely define the corresponding MOs.
This question is treated in \cite{KoSe} and \cite{KoPuMaDe}.

\begin{defi}\label{deficanon}
  {\rm Given a SP of length $d+1$ we define the {\em canonical} MO
    corresponding to it as follows. The SP is read from the back and to each
    encountered couple of equal (resp. different) consecutive signs one puts
    in correspondence the letter $N$
    (resp. $P$) in the MO. E.g. for $d=7$, to the SP $(+,+,-,-,+,-,+,+,-)$
    there corresponds the canonical MO $P<N<P<P<P<N<P<N$. The canonical MO
    is obtained when one constructs a HP realizing the given SP using consecutive 
    products of the form $(x\pm \varepsilon )Q(x)$,
    see Remark~\ref{remconcat}. Each SP is realizable by its canonical MO, see
    \cite{KoSe}. A SP is called {\em canonical} if it is realizable only by
  its canonical~MO.}
\end{defi}

For SPs one can use the notation $\Sigma _{p_1,p_2,\ldots ,p_s}$, where $p_i$
are the lengths of the maximal sequences of equal signs. E.g. the SP in
Definition~\ref{deficanon} is $\Sigma _{2,2,1,1,2,1}$. The following necessary
condition for a SP to be canonical is proved in~\cite{KoSe}:

\begin{tm}
  If the SP $\Sigma _{p_1,p_2,\ldots ,p_s}$ is canonical, then there are no two
  consecutive numbers $p_i$ which are larger than $1$, and for
  $2\leq i\leq s-1$, one has $p_i\neq 2$.
\end{tm}

\begin{rems}
  {\rm (1) Thus for $d\geq 3$, the SPs $\Sigma _{\pm}$ (corresponding to rigid MOs,
see Notation~\ref{nota1} and Theorem~\ref{tmmain}) are not canonical. For $d=1$ and $2$, they are canonical, see Example~\ref{ex1}.

(2) The SPs with $c=d$, $p=0$ and $c=0$, $p=d$, are canonical.
They correspond to the trivial case when all roots are positive or negative, see the lines after Definition~\ref{defirigid}..

(3) The MO
corresponding to a canonical SP for which one does not have $s=1$ or
$p_1=\cdots =p_s=1$, 
with at least one number $p_i$ larger than $2$ (or with $p_1=2$
or with $p_s=2$), is not rigid. Indeed, the
presence of a number $p_i>2$ for $2\leq i\leq s-1$
(or of $p_1>1$ or of $p_2>1$) implies the presence of $p_i-1\geq 2$
(or of $p_1$ or of $p_s$) consecutive
letters $P$ or $N$ in the MO, see Definition~\ref{deficanon}.
Thus in and only in the trivial case does one have a rigid MO realizing a canonical~SP.

(4) The SPs of the form $\Sigma _{1,p_2}$, $\Sigma _{p_1,1}$,
$\Sigma _{1,p_2,1}$, $p_2\geq 3$, and $\Sigma _{p_1,1,p_3}$ are canonical,
see~\cite{KoPuMaDe}.}
  \end{rems}

The problems treated in the present paper are part of problems about
real (not necessarily hyperbolic) univariate polynomials. For such a polynomial without
vanishing coefficients, Descartes' rule of signs implies that
the number $pos$ of its positive roots is not greater
than the number $c$ of sign changes in the sequence of its coefficients, and
the difference $c-pos$ is even. In the same way, for the number $neg$ of its
negative roots, one has $neg\leq p$ and $p-neg\in 2\mathbb{Z}$, where $p$
is the number of sign preservations.

The problem for which couples $(pos, neg)$
compatible with these requirements can one find such a real polynomial
with prescribed signs of its coefficients seems to have been formulated
for the first time in~\cite{AJS}. For $d=4$, D.~Grabiner has obtained the
first nontrivial result, i.e. a compatible, but not realizable couple
$(pos, neg)$, see~\cite{Gr}. In the cases $d=5$ and $6$ the problem has been
thoroughly studied in \cite{AlFu} while the exhaustive answer for $d=7$ and $8$
can be found in \cite{FoKoSh} and \cite{KoCzMJ}. For $d\leq 8$, all
compatible, but not realizable cases, are ones in which either $pos=0$ or
$neg=0$.

For $d\geq 9$, there are examples of compatible and nonrealizable couples
$(pos, neg)$ with $pos\geq 1$ and $neg\geq 1$, see \cite{KoMB} and~\cite{CGK}.
Various problems about HPs are exposed in~\cite{Ko}. A tropical
  analog of Descartes' rule of signs is discussed in~\cite{FoNoSh}.

\section{Proof of Theorem~\protect\ref{tmmain}\protect\label{secprtmmain}}

\begin{proof}[Proof of part (1)]
  Suppose that for $d\geq 3$, a MO $r$ contains the string of
  inequalities $P<P<N$.
  Consider the two polynomials

  $$\begin{array}{llllll}
    P_1&:=&(x-1)(x-1.1)(x+3)&=&x^3+0.9x^2-5.2x+3.3&{\rm and}\\ \\
    P_2&:=&(x-1)(x-3)(x+3.1)&=&x^3-0.9x^2-9.4x+9.3~.\end{array}$$
  They define two different SPs: $\sigma (P_1)=(+,+,-,+)$ and
  $\sigma (P_2)=(+,-,-,+)$. Hence one can
  realize the whole MO $r$ by two different SPs starting with the
  polynomials $P_1$ and $P_2$ and using $d-3$ multiplications with one and the same
  polynomials $x\pm \varepsilon$ or $1\pm \varepsilon x$, see
  Remark~\ref{remconcat}. After each multiplication one obtains again two polynomials defining different SPs.  Hence $r$ is not rigid.

  If the MO contains a string of inequalities $N<N<P$, $N<P<P$ or
  $P<N<N$, then one can consider instead of the polynomials $P_j$, $j=1$, $2$,
  the polynomials $S_j:=-P_j(-x)$, $T_j:=x^3P_j(1/x)$ and $R_j:=x^3P_j(-1/x)$
  respectively and perform a similar reasoning. The SPs defined by these
  polynomials are:

  $$\begin{array}{lll}
    \sigma (S_1)=(+,-,-,-)~,&\sigma (S_2)=(+,+,-,-)~,&
    \sigma (T_1)=(+,-,+,+)~,\\ \\
    \sigma (T_2)=(+,-,-,+)~,&\sigma (R_1)=(+,+,+,-)~,&\sigma (R_2)=(+,+,-,-)~,
    \end{array}$$
 hence $\sigma (S_1)\neq \sigma (S_2)$,  $\sigma (T_1)\neq \sigma (T_2)$ and $\sigma (R_1)\neq \sigma (R_2)$.
  \end{proof}

\begin{proof}[Proof of part (2)]
  We prove part (2) of the theorem by induction on $d$.
  For $d=1$ and $2$, the theorem is to be checked straightforwardly, see
  Example~\ref{ex1}. Suppose
  that part (2) of the theorem holds true for $d\leq d_0$, $d_0\geq 2$. Set $d:=d_0+1$.
  The sign of the
  constant term of a HP realizing the given MO depends only on the
  signs of the roots, not on the MO. So this sign is also to be checked
  directly.

  Consider a polynomial $Q:=x^{d_0+1}+\sum _{j=0}^{d_0}b_jx^j$
  defining the given MO $\rho$ with $d=d_0+1$, with $\rho$ standing for
  $r_{PP}$, $r_{PN}$, $r_{NP}$ or $r_{NN}$. We represent
  it in the form

  $$Q:=(x-\varphi )(x-\psi )V~,~~~{\rm where}~~~V:=\prod _{j=1}^{d_0-1}(x-\xi _j)=
  x^{d_0-1}+\sum _{j=0}^{d_0-2}c_jx^j~.$$
  Here $\varphi$ and $\psi$ are the two roots of $Q$ of least moduli,
  $|\varphi |<|\psi |$, and $\xi _j$ are its other roots. The signs of
  $\varphi$ and $\psi$ are opposite. Denote by $r$ the MO defined by the
  polynomial $V$. Using the notation of Remark~\ref{remconcat} one can
  say that the MO defined by the polynomial $R:=(x-\psi )V$ is either
  $Pr$ or $Nr$ depending to the sign of the root $\psi$, and the MO $\rho$
  is either $NPr$ or $PNr$.

  We denote by $\Sigma$ the SP $\Sigma _+$ or $\Sigma _-$ according to the case
  and by $\Sigma '$ and $\Sigma ''$ the SPs obtained from $\Sigma$ by deleting
  its one or two last signs respectively.

  We include $Q$ into a one-parameter family of polynomials of the form

  $$Z_t:=(x+t\psi )(x-\psi )\prod _{j=1}^{d_0-1}(x-\xi _j)~~~,~~~t\in [0,1]~.$$
  As $\varphi \cdot \psi <0$ and $|\varphi |<|\psi |$,
  there exists $t_*\in (0,1)$ such that $\varphi =-t_*\psi$, i.e. $Z_{t_*}=Q$.

  For $t=0$, one obtains $Z_0=xR$. The theorem being true for $d=d_0$ and
  $d=d_0-1$, the
  polynomial $R$ defines the SP $\Sigma '$,
  because $R$ defines the MO $Pr$ or~$Nr$, and $V$ defines the
  MO $r$ and the SP~$\Sigma ''$.

  For $t=1$, one has $Q=(x^2-\psi ^2)V$. Hence

  $$Z_1=x^{d_0+1}+c_{d_0-2}x^{d_0}+\left( \sum _{j=0}^{d_0}(c_j-\psi ^2c_{j+2})x^{j+2}\right) -
  \psi ^2(c_1x+c_0)~.$$
  The signs of $c_j$ and $c_{j+2}$ are opposite (see Notation~\ref{nota1} for
  the definition of the SPs
  $\Sigma _{\pm}$), therefore sgn$(c_j-\psi ^2c_{j+2})=$sgn$(c_j)$. Thus the first
  $d_0$ coefficients of $Z_1$ have the signs given by the SP $\Sigma$.
  This is the case of the last two coefficients as well, because
  sgn$(-\psi ^2c_1)=-$sgn$(c_1)$ and sgn$(-\psi ^2c_0)=-$sgn$(c_0)$. Hence
  $Z_1$ defines the SP~$\Sigma$.

  The coefficients of $Z_t$ are linear functions in $t\in [0,1]$. If
  their signs for $t=0$ and $t=1$ are the corresponding components of the
  SP $\Sigma$, then this is the case for any $t\in [0,1]$. (For the constant
  term, one has to consider its values for $t=1$
  and for $t>0$ close to $0$.) In particular, for $t=t_*$, the signs are the ones of
  the SP~$\Sigma$. This proves part~(2) of the theorem.
  
  \end{proof}

\section{Proof of Theorem~\protect\ref{tm2}\protect\label{secprtm2}}

Without loss of generality we limit ourselves to the case of monic HPs. 
We prove the theorem by induction on $d$. The cases $d=1$ and $2$ are
considered in Examples~\ref{ex1} and \ref{ex2}. For $d=1$, no 
coefficient of the HP equals~$0$.

Suppose now that $d\geq 3$. 
We assume that there is at least one equality
between a modulus of a negative
root and a positive root, otherwise one can apply Theorem~\ref{tmmain}. So
suppose that the HP has roots $\pm a$, $a\neq 0$, and the HP
is of the form $S:=(x^2-a^2)Q$, where $Q$ is a degree $d-2$
HP without
root at~$0$. Thus the roots of $Q$ define one of the MOsAE
$r_{PN}^0$, $r_{PP}^0$, $r_{NN}^0$ and $r_{NP}^0$, so one can use
the inductive assumption.

If $d$ is odd, then $Q$ has no vanishing coefficient
and defines one of the SPs $\Sigma _{\pm}$. Set $Q:=\sum _{j=0}^{d-2}q_jx^j$. Then

\begin{equation}\label{eqS}
  S=q_{d-2}x^d+q_{d-3}x^{d-1}+\left( \sum _{j=0}^{d-4}(q_j-a^2q_{j+2})x^{j+2}
  \right) -a^2q_1x-a^2q_0~.
\end{equation}
The first two and the last two of the coefficients of $S$ are obviously
nonzero. For the others one can observe that as by inductive assumption
$Q$ defines one of the SPs $\Sigma _{\pm}$ hence
$q_j\cdot q_{j+2}<0$, one has $q_j-a^2q_{j+2}\neq 0$. This proves part (1)
of the theorem.

If $d$ is even, then $Q$ can have a vanishing coefficient in which case $Q$ is
of the form $\prod _{j=1}^{(d-2)/2}(x^2-a_j^2)$ hence $S$ is of the form
$\prod _{j=1}^{d/2}(x^2-a_j^2)$ (with $a_{d/2}=a$) and defines the SPAZ
  $(+,0,-,0,+,0,-,0,\ldots )$.

If $d$ is even and $Q$ has no vanishing
coefficient, then the set of its roots is not representable as a union of
couples $\{ a_j, -a_j\}$, $a_j\in \mathbb{R}^*$, so this is the case of $S$
as well. Moreover, using equality (\ref{eqS}) in the same way as for $d$ odd
one concludes that $S$ has no vanishing coefficient.
Part (2) of the theorem is proved.

\end{document}